\documentclass{pnastwo}

\usepackage{amssymb,amsfonts,amsmath}

\usepackage{amssymb}
\usepackage{amstext}
\usepackage{amsmath}
\usepackage{amscd}
\usepackage{amsfonts}
\usepackage{enumerate}
\usepackage{graphicx}
\usepackage{latexsym}
\usepackage{mathrsfs}
\usepackage[all,color,dvips]{xy}
\usepackage{mathtools}
\xyoption{all}

\usepackage{pstricks}
\usepackage{lscape}
\usepackage{comment}

\newcommand{\DDD}{\mathsf{D}}
\newcommand{\bo}{\operatorname{b}\nolimits}
\newcommand{\gl}{\operatorname{gl.dim}\nolimits}
\DeclareMathOperator{\CM}{\mathsf{CM}}

\newcommand{\Hom}{\operatorname{Hom}\nolimits}
\renewcommand{\mod}{\mathsf{mod}\hspace{.01in}}
\newcommand{\proj}{\mathsf{proj}\hspace{.01in}}
\newcommand{\inj}{\mathsf{inj}\hspace{.01in}}

\newcommand{\xto}{\xrightarrow}

\newcommand{\End}{\operatorname{End}\nolimits}
\newcommand{\Ext}{\operatorname{Ext}\nolimits}
\newcommand{\op}{\operatorname{op}\nolimits}
\newcommand{\ann}{\operatorname{ann}\nolimits}

\newcommand{\Tr}{\operatorname{Tr}\nolimits}

\newcommand{\CC}{{\mathcal C}}
\newcommand{\TT}{{\mathcal T}}
\newcommand{\FF}{{\mathcal F}}
\newcommand{\XX}{{\mathcal X}}
\newcommand{\YY}{{\mathcal Y}}
\newcommand{\EE}{{\mathcal E}}

\newcommand{\add}{\mathsf{add}\hspace{.01in}}
\newcommand{\Fac}{\mathsf{Fac}\hspace{.01in}}
\newcommand{\Sub}{\mathsf{Sub}\hspace{.01in}}
\newcommand{\KKb}{\mathsf{K}^{\rm b}}

\newcommand{\thick}{\mathsf{thick}\hspace{.01in}}

\newcommand{\ftors}{\mbox{\rm f-tors}\hspace{.01in}}

\newcommand{\fftors}{\mbox{\rm ff-tors}\hspace{.01in}}

\newcommand{\twosilt}{\mbox{\rm 2-silt}\hspace{.01in}}

\newcommand{\ctilt}{\mbox{\rm c-tilt}\hspace{.01in}}
\newcommand{\rigid}{\mbox{\rm rigid}\hspace{.01in}}
\newcommand{\mrigid}{\mbox{\rm m-rigid}\hspace{.01in}}
\newcommand{\trigid}{\mbox{\rm $\tau$-rigid}\hspace{.01in}}
\newcommand{\sttilt}{\mbox{\rm s$\tau$-tilt}\hspace{.01in}}

\newcommand{\ttilt}{\mbox{\rm $\tau$-tilt}\hspace{.01in}}

\newcommand{\tilt}{\mbox{\rm tilt}\hspace{.01in}}

\contributor{Submitted to Proceedings
of the National Academy of Sciences of the United States of America}
\url{www.pnas.org/cgi/doi/10.1073/pnas.0709640104}
\copyrightyear{2008}
\issuedate{Issue Date}
\volume{Volume}
\issuenumber{Issue Number}

\begin{document}



\title{Introduction to $\tau$-tilting theory}





\author{Osamu Iyama\affil{1}{Nagoya University, Nagoya, Japan},
Idun Reiten\affil{2}{NTNU, Trondheim, Norway}}

\contributor{Submitted to Proceedings of the National Academy of Sciences
of the United States of America}

\maketitle

\begin{article}

\begin{abstract}
From the viewpoint of mutation, we will give a brief survey of tilting theory and cluster-tilting theory together with a motivation from cluster algebras.
Then we will give an introdution to $\tau$-tilting theory which was recently developed in \cite{AIR}.
\end{abstract}







\section{1. Introduction}

Let $\Lambda$ be a finite dimensional algebra over an algebraically closed field $k$,
for example $k$ is the field of complex numbers.
We always assume that $\Lambda$ is associative and has an identity.
An important class of algebras is the \emph{path algebras} $kQ$ for
a finite \emph{quiver} $Q$, that is, a finite directed graph.
The $k$-basis for $kQ$ is the paths in $Q$.
For example, when $Q$ is $1\xrightarrow{a}2\xrightarrow{b}3$,
then the paths are $\{e_1,e_2,e_3,a,b,ba\}$. Here $e_1,e_2,e_3$
are the trivial paths, starting and ending at 1,2,3 respectively.
Then multiplication in $kQ$ is induced by composition of paths, so that
for example $b\cdot a=ba$, $a\cdot b=0$, $a\cdot e_1=a$, $a\cdot e_2=0$.
When $Q$ is \emph{acyclic} (i.e. contains no oriented cycles), then
$kQ$ is a finite dimensional $k$-algebra.

The representation theory of $\Lambda$ deals with the investigation of the module 
category $\mod\Lambda$ of finitely generated left $\Lambda$-modules. 
We refer to \cite{ARS,ASS,Rin1} for background on representation theory.
One of the basic properties of $\mod\Lambda$ is the \emph{Krull-Schmidt} property:
Any module $X$ is isomorphic to a finite direct sum $X_1\oplus\cdots\oplus X_m$ of
indecomposable modules $X_1,\ldots,X_m$, which are uniquely determined
up to isomorphism and permutation.
We say that $X$ is \emph{basic} if $X_1,\ldots,X_m$ are pairwise nonisomorphic,
and we write $|X|=m$ in this case.
We assume that $\Lambda$ is basic as a $\Lambda$-module.
Then $\Lambda$ is isomorphic to $kQ/I$ for a finite quiver $Q$ and an ideal $I$ in $kQ$.
In this case, let $\{1,\ldots,n\}$ be the set of vertices in $Q$. Then
$|\Lambda|=n$ holds, and we have
a decomposition $\Lambda=P_1\oplus\cdots\oplus P_n$, where $P_i=\Lambda e_i$ is
an indecomposable projective $\Lambda$-module. Moreover the simple $\Lambda$-modules are
$S_1,\ldots,S_n$, where $S_i=P_i/{\rm rad} P_i$ for the radical ${\rm rad} P_i$ of $P_i$.
A famous theorem of Gabriel tells us that a path algebra $kQ$ of
a connected finite acyclic quiver $Q$ has only a finite number of nonisomorphic indecomposable modules if and only if the underlying 
graph of $Q$ is one of the simply laced Dynkin diagrams $A_n$, $D_n$, $E_6$,
$E_7$, $E_8$.

An important class of modules has been the \emph{tilting modules},
introduced more than 30 years ago \cite{BB,HR,B} as a simultaneous
generalization of progenerators in classical Morita theory and
BGP-reflections for quiver representations \cite{BGP}. They are a certain class of
basic modules $T=T_1\oplus\cdots\oplus T_n$ with $n=|\Lambda|$,
where each $T_i$ is indecomposable. BGP-reflections can be realized by a special
class of tilting modules called APR-tilting modules \cite{APR}.
One of the important properties of tilting modules is that when we remove some
direct summand $T_i$ from $T$ to get $T/T_i=\bigoplus_{j\ne i}T_j$ (which is called an
\emph{almost complete tilting module}),
then there is at most one indecomposable module $T_i^*$
such that $T_i^*\not\simeq T_i$ and $(T/T_i)\oplus T_i^*$ is again a tilting module
(which is called \emph{mutation} of $T$).
Then $T_i$ (and $T_i^*$ if it exists) are said to be the \emph{complement(s)} of $T/T_i$.
Thus any almost complete tilting module has either one or two complements
(Theorem \ref{complements for tilting}), and mutation is possible only when
there are two complements.

To make mutation always possible, it is desirable to enlarge our class of
tilting modules in order to get the more regular property that
almost complete ones always have two complements.
This is accomplished by introducing \emph{$\tau$-tilting modules},
or more precisely, \emph{support $\tau$-tilting pairs} \cite{AIR}.
Such a property is of central interest when dealing with
categorification of cluster algebras. For cluster categories \cite{BMRRT},
and more generally, for 2-Calabi-Yau triangulated categories \cite{IY},
a class of objects called \emph{cluster-tilting objects} have the desired property
(Theorem \ref{complements for cluster-tilting}).
For path algebras $kQ$, this result for cluster categories was
translated to the fact that the same holds for \emph{support tilting
pairs} (which generalize tilting modules) \cite{IT,Rin3}
(Theorem \ref{support tilting}). More
generally, for 2-Calabi-Yau triangulated categories we can translate to the
property that for the associated 2-Calabi-Yau tilted algebras
the support $\tau$-tilting pairs have the desired property for
complements (see section 5). This raised the question of a possible extension
to support $\tau$-tilting pairs for other classes of algebras.
In fact it turned out to be the case for any finite dimensional algebra
\cite{AIR} (Theorem \ref{2 complements}).

For any finite dimensional $k$-algebra $\Lambda$, in particular for the path algebras
of Dynkin quivers, we have an associated AR-quiver, where the vertices correspond to the
isomorphism classes of the indecomposable modules. The arrows correspond to
what is called irreducible maps, and the dotted lines indicate
a certain functor $\tau$ called AR-translation (see section 4.1 for details).
For example, for the path algebra $kQ$ of the quiver
$1\xrightarrow{a}2\xrightarrow{b}3$, we have the following AR-quiver
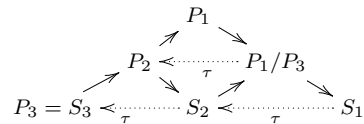
\[{\small\xymatrix@R=.5em@C=.8em{
&&P_1\ar[dr]\\
&P_2\ar[ur]\ar[dr]&&P_1/P_3\ar@{.>}[ll]^\tau\ar[dr]\\
P_3=S_3\ar[ur]&&S_2\ar@{.>}[ll]^\tau\ar[ur]&&S_1\ar@{.>}[ll]^\tau
}}\]
In this case there are 5 tilting modules
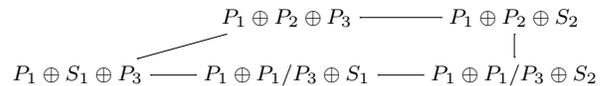
\[{\small\xymatrix@R1em@C2em{
&P_1\oplus P_2\oplus P_3\ar@{-}[r]&P_1\oplus P_2\oplus S_2\ar@{-}[d]\\
P_1\oplus S_1\oplus P_3\ar@{-}[ru]\ar@{-}[r]&P_1\oplus P_1/P_3\oplus S_1\ar@{-}[r]&
P_1\oplus P_1/P_3\oplus S_2
}}\]
where the edges indicate mutation. This picture is incomplete since, for example,
it is impossible to replace $P_1$ in the tilting module $P_1\oplus P_2\oplus P_3$
to get another tilting module.
Consider now a slightly bigger category called the \emph{cluster category}
(see section 3) with the following AR-quiver:
{\small
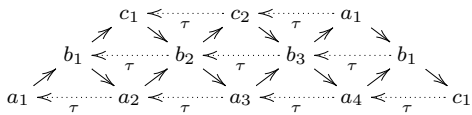
\[
\xymatrix@R=.5em@C=.8em{
&&c_1\ar[dr]&&c_2\ar[dr]\ar@{.>}[ll]^\tau&&a_1\ar[dr]\ar@{.>}[ll]^\tau\\
&b_1\ar[ur]\ar[dr]&&b_2\ar[ur]\ar[dr]\ar@{.>}[ll]^\tau&&
b_3\ar[ur]\ar[dr]\ar@{.>}[ll]^\tau&&b_1\ar[dr]\ar@{.>}[ll]^\tau\\
a_1\ar[ur]&&a_2\ar[ur]\ar@{.>}[ll]^\tau&&a_3\ar[ur]\ar@{.>}[ll]^\tau&&
a_4\ar[ur]\ar@{.>}[ll]^\tau&&c_1\ar@{.>}[ll]^\tau
}\]}
In this category there is an important class of objects called
\emph{cluster-tilting objects}, which generalize tilting modules.
In this case there are 14 cluster-tilting objects
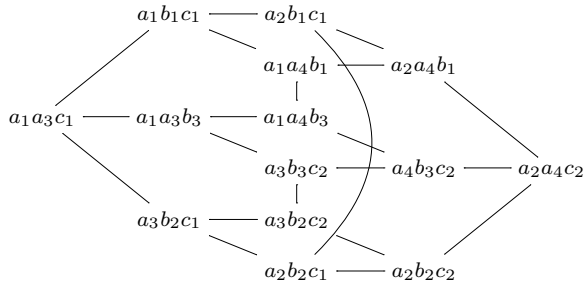
\[{\small\xymatrix@R.8em@C2em{
&a_1b_1c_1\ar@{-}[r]\ar@{-}[rd]&a_2b_1c_1\ar@{-}[rd]\\
&&a_1a_4b_1\ar@{-}[r]\ar@{-}[d]&a_2a_4b_1\ar@{-}[rdd]\\
a_1a_3c_1\ar@{-}[ruu]\ar@{-}[r]\ar@{-}[rdd]&a_1a_3b_3\ar@{-}[r]\ar@{-}[rd]&
a_1a_4b_3\ar@{-}[rd]\\
&&a_3b_3c_2\ar@{-}[r]\ar@{-}[d]&a_4b_3c_2\ar@{-}[r]&a_2a_4c_2\\
&a_3b_2c_1\ar@{-}[r]\ar@{-}[rd]&
a_3b_2c_2\ar@{-}[rd]\\
&&a_2b_2c_1\ar@{-}[r]\ar@{-}@/_10mm/[uuuuu]&a_2b_2c_2\ar@{-}[ruu]
}}\]
where the symbol $\oplus$ is omitted and the edges indicate mutation.
This picture is complete in the sense that
we can replace any element in any cluster-tilting object to get another
cluster-tilting object.
This corresponds to the fact that we have exactly two complements for
almost complete cluster-tilting objects in cluster categories
(see Theorem \ref{complements for cluster-tilting}).

\section{2. Tilting theory}

In this section we include some main results on tilting modules and
related objects, focusing on when we have a well defined operation
of mutation on them.
We refer to \cite{ASS,Rin1,Ha,AHK} for more background on tilting theory.

\subsection{2.1. Properties of tilting modules}

\ Let as before $\Lambda$ be a finite dimensional $k$-algebra with $n=|\Lambda|$.
A $\Lambda$-module $T$ is a \emph{partial tilting module} if
$T$ has projective dimension at most one and $\Ext^1_\Lambda(T,T)=0$.
It is a \emph{tilting module} if moreover $|T|=n$. This condition can be
replaced by the condition: There is an exact sequence
$0\to\Lambda\to T_0\to T_1\to0$ with $T_0$ and $T_1$ in $\add T$.
Here $\add T$ denotes the subcategory of $\mod\Lambda$ consisting of
direct summands of finite direct sums of copies of $T$.
The following result due to Bongartz is basic in tilting theory.

\begin{theorem}\emph{\cite{B}}\label{Bongartz}
Any partial tilting module is a direct summand of a tilting module.
\end{theorem}

It is natural to ask how many tilting modules exist for a given partial tilting module.
There is an explicit answer for an \emph{almost complete tilting module}, that is, 
a partial tilting module $U$ satisfying $|U|=n-1$.
An indecomposable module $X$ is a \emph{complement} for an almost complete tilting module
$U$ if $U\oplus X$ is a tilting module.
The following result is basic in this paper, where a $\Lambda$-module $X$
is called \emph{faithful} if the annihilator
$\ann X:=\{a\in\Lambda\mid aX=0\}$ of $X$ is zero.

\begin{theorem}\label{complements for tilting}\emph{\cite{HU1,RS,U}}
Let $U$ be an almost complete tilting module.\\
{\rm (a)} $U$ has either one or two complements.\\
{\rm (b)} $U$ has two complements if and only if $U$ is faithful.
\end{theorem}

When there are exactly two complements, then there is a nice relationship between them.

\begin{theorem}\label{exchange sequence for tilting}\emph{\cite{HU1,RS}}
Let $X$ and $Y$ be two complements for an almost complete tilting
module $U$. After we interchange $X$ and $Y$ if necessary,
there is an exact sequence
$0\to X\xrightarrow{f}U'\xrightarrow{g}Y\to0$ where $f$ is a
minimal left $(\add U)$-approximation and $g$ is a
minimal right $(\add U)$-approximation.
\end{theorem}

Recall that $g:U'\to Y$ is a \emph{right $(\add U)$-approximation}
if $U'\in\add U$ and any map $h:U\to Y$ factors through $g$.
Moreover it is called (right) \emph{minimal} if any map $a:U'\to U'$
satisfying $ga=g$ is an isomorphism.

Let again $Q$ be the quiver $1\to 2\to 3$ and
$kQ$ the corresponding path algebra. 
The vertex 3 of the quiver $Q$ is a \emph{sink}, that is, no arrow starts
at 3. Then a \emph{reflection} (or \emph{mutation}) of $Q$ at vertex 3 is
defined by reversing all arrows ending at 3, to get the quiver
$Q'=\mu_3(Q): 1\xrightarrow{}2\leftarrow{}3$ in this example \cite{BGP}.
We have the following close relationship between $\mod kQ$ and $\mod kQ'$.

\begin{theorem}\label{BGP thm}\emph{\cite{BGP}}
Let $Q$ be a finite acyclic quiver with a sink $i$ and $Q':=\mu_i(Q)$. Then
there is an equivalence of categories $F:\EE_{S_i}\to\EE_{S'_i}$, where
$\EE_{S_i}$ (respectively, $\EE_{S'_i}$) is the full subcategory of
$\mod kQ$ (respectively, $\mod kQ'$) consisting of modules without direct summands $S_i$ (respectively, $S'_i$).
\end{theorem}

It was shown in \cite{APR} that the $kQ$-module $T=(kQ/S_{i})\oplus\tau^{-1}S_{i}$ has the endomorphism algebra $kQ'$, and the above functor $F$ is isomorphic to $\Hom_{kQ}(T,-)$.
This $T$ is a special case of what was later defined to be tilting modules.

An important notion in tilting theory is a \emph{torsion pair}, that is, a pair $(\TT,\FF)$ of full subcategories of $\mod\Lambda$ such that (i) $\Hom_\Lambda(\TT,\FF)=0$ and (ii) for any $X$ in $\mod\Lambda$, there exists an exact sequence $0\to Y\to X\to Z\to 0$ with $Y\in\TT$ and $Z\in\FF$.
In this case $\TT$ is a \emph{torsion class}, that is, a subcategory closed under
extensions, factor modules and isomorphisms, and $\FF$ is a \emph{torsionfree class},
that is, a subcategory closed under extensions, submodules and isomorphisms.
For any $M$ in $\mod\Lambda$, we have a torsion class
${}^\perp M:=\{X\in\mod\Lambda\mid\Hom_\Lambda(X,M)=0\}$ and a torsionfree class
$M^\perp:=\{X\in\mod\Lambda\mid\Hom_\Lambda(M,X)=0\}$.
Let $\Fac M$ (respectively, $\Sub M$) be the subcategory of $\mod\Lambda$ consisting of factor modules (respectively, submodules) of direct sums of copies of $M$.
The equivalence in Theorem \ref{BGP thm} has the following interpretation in tilting theory:

\begin{theorem}\label{tilting and torsion}\emph{\cite{BB,HR}}
Let $T$ be a tilting module over a finite dimensional $k$-algebra $\Lambda$ and $\Gamma:=\End_\Lambda(T)^{\op}$.\\
{\rm (a)} $\TT:=\Fac T$ and $\FF:=T^\perp$ give a torsion pair $(\TT,\FF)$ in $\mod\Lambda$,
and $\XX:={}^\perp DT$ and $\YY:=\Sub DT$ (where $D$ is the $k$-dual) give a torsion pair
$(\XX,\YY)$ in $\mod\Gamma$.\\
{\rm (b)} We have equivalences of categories $\Hom_\Lambda(T,-):\TT\to\YY$
and $\Ext^1_\Lambda(T,-):\FF\to\XX$.
\end{theorem}

Note that in the above special case $S_3$ is the only indecomposable $kQ$-module
not in $\TT=\Fac T$, and $\FF=\add S_3$.

There is a close connection between tilting modules and torsion classes.
We need an important notion of \emph{functorially finite subcategories} \cite{AS}.
For torsion classes $\TT$, being functorially finite is equivalent to the existence
of $M$ in $\mod\Lambda$ such that $\TT=\Fac M$.
We denote by $\tilt\Lambda$ the set of isomorphism classes of basic tilting $\Lambda$-modules,
and by $\fftors\Lambda$ the set of faithful functorially finite torsion classes in $\mod\Lambda$.

\begin{theorem}\label{tilting and torsion2}\emph{\cite{AS,Ho,S,As}}
Let $\Lambda$ be a finite dimensional $k$-algebra. Then there is a bijection $\tilt\Lambda\to\fftors\Lambda$ given by $T\mapsto\Fac T$.
\end{theorem}

Let $(T,P)$ be a pair of $\Lambda$-modules, where $P$ is projective,
$\Hom_\Lambda(P,T)=0$ and $T$ is a partial tilting $\Lambda$-module.
Then $(T,P)$ is called a \emph{support partial tilting pair} for $\Lambda$.
It is called a \emph{support tilting pair} (respectively, \emph{almost complete support
tilting pair}) if $|T|+|P|=n$ (respectively, $|T|+|P|=n-1$), where $n=|\Lambda|$. 
For an indecomposable module $X$, we say that $(X,0)$ (respectively, $(0,X)$)
is a \emph{complement} for an almost complete support tilting pair $(T,P)$ if
$(T\oplus X,P)$ (respectively, $(T,P\oplus X)$) is a support tilting pair.
We then have the following result for path algebras.

\begin{theorem}\emph{\cite{IT,Rin3}}\label{support tilting}
Let $\Lambda=kQ$ be a path algebra of an acyclic quiver $Q$. Then any almost complete
support tilting pair has precisely two complements.
\end{theorem}

Note that Theorem \ref{support tilting} is not true for finite dimensional $k$-algebras
in general. For example,
let $Q$ be the quiver $1\xrightarrow{a}2\xrightarrow{b}3$ and $\Lambda:=kQ/\langle ba\rangle$.
In this case there are 8 support tilting pairs for $\Lambda$
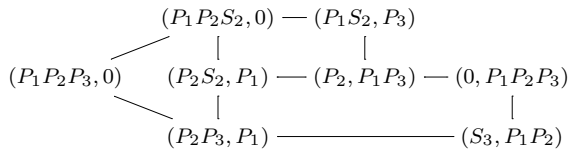
\[{\small\xymatrix@R1em@C1em{
&(P_1P_2S_2,0)\ar@{-}[r]\ar@{-}[d]&(P_1S_2,P_3)\ar@{-}[d]&\\
(P_1P_2P_3,0)\ar@{-}[ru]\ar@{-}[rd]&(P_2S_2,P_1)\ar@{-}[r]&
(P_2,P_1P_3)\ar@{-}[r]&(0,P_1P_2P_3)\\
&(P_2P_3,P_1)\ar@{-}[u]\ar@{-}[rr]&&(S_3,P_1P_2)\ar@{-}[u]
}}\]
where the symbol $\oplus$ is omitted and edges indicate mutation.
For example, the almost complete support tilting pair
$(P_1\oplus P_3,0)$ does not have two complements.

There are two important quivers whose vertices are tilting modules.
One is the \emph{Hasse quiver} of $\tilt\Lambda$, where we regard $\tilt\Lambda$ as a partially ordered set by defining $T\ge U$ when $\Fac T\supset\Fac U$.
Thus we draw an arrow $T\to U$ if $T>U$ and there is no $V\in\tilt\Lambda$ satisfying $T>V>U$. 
The other is the \emph{exchange quiver} of $\tilt\Lambda$. We draw an arrow $X\oplus U\to Y\oplus U$ when $X$ and $Y$ are complements of $U$ and there is an exact sequence $0\to X\to U'\to Y\to0$ as given in Theorem \ref{exchange sequence for tilting}.
Here we have the following.

\begin{theorem}\emph{\cite{HU2}}\label{Hasse and exchange}
The Hasse quiver and the exchange quiver coincide.
\end{theorem}

\subsection{2.2. Generalizations of tilting modules}

\ Tilting modules discussed in the previous subsection are often called `classical'
tilting modules, and there are many generalizations.
Here we are mainly interested in cases where there is some kind
of mutation, which is given via approximation sequences.

We have a more general class of `tilting modules of finite projective dimension', which
by definition are $\Lambda$-modules $T$ with the following properties:
(i) $T$ has finite projective dimension, (ii) $\Ext^j_\Lambda(T,T)=0$ for $j>0$,
(iii) there is an exact sequence $0\to\Lambda\to T_0\to
\cdots\to T_m\to0$, where all $T_i$ are in $\add T$ \cite{Miy,Ha}.
The notion of almost complete tilting modules and their complements are defined
in the same way as before.
In this case there are in general more than two complements, and all the complements are
usually connected by approximation sequences \cite{CHU}.
There are other aspects of such tilting modules which have been
investigated more, like connections with functorially finite subcategories,
with applications to algebraic
groups and Cohen-Macaulay representations \cite{AB,AR,Rin2}.

One of the important properties of tilting modules of finite projective dimension is that
they induce a derived equivalence between $\Lambda$ and $\End_\Lambda(T)^{\op}$.
From this viewpoint of derived Morita theory, the notion of tilting complexes introduced
in \cite{Ri} is most important. Also here there are usually more than two complements
for almost complete ones.
From the viewpoint of mutation theory, the notion of silting complexes \cite{KV}
is more natural than that of tilting complexes.
We denote by $\KKb(\proj\Lambda)$ the homotopy category
of bounded complexes of finitely generated projective modules.
A \emph{tilting complex} (respectively, \emph{silting complex})
is an object $P$ in $\KKb(\proj\Lambda)$ satisfying
$\Hom_{\KKb(\proj\Lambda)}(P,P[i])=0$ for any $i\neq0$ (respectively, $i>0$)
and $\thick P=\KKb(\proj\Lambda)$. Here $\thick P$ is the smallest
full subcategory of $\KKb(\proj\Lambda)$ which contains $P$ and
is closed under cones, $[\pm1]$, direct summands and isomorphisms.
There is mutation theory analogous to the case of tilting modules \cite{AI}. 
One has operations of left and right mutations obtained using approximation sequences,
which shows that there are always infinitely many complements for an almost complete silting complex.
Moreover there is a natural partial order whose Hasse quiver coincides with the exchange quiver.
If we restrict to silting complexes which are `two-term', then almost complete
silting complexes have exactly two complements \cite{AIR}
(Corollary \ref{complement for two-term silting}).

\section{3. Cluster algebras and 2-Calabi-Yau categories}

In this section we discuss cluster categories and 2-Calabi-Yau triangulated categories,
which were motivated by trying to categorify the essential ingredients of the cluster algebras
of Fomin-Zelevinsky \cite{FZ1,FZ2}.
We restrict to considering what is called acyclic cluster algebras, with no coefficients.
They are associated with a finite acyclic quiver with $n$ vertices.
Essential concepts in the theory of cluster algebras are clusters, cluster variables,
seeds, and quiver mutation.

Let $Q$ be a finite connected quiver having no loops or 2-cycles. For a vertex $i$ of $Q$,
we define a new quiver $\mu_i(Q)$ having no loops or 2-cycles, called \emph{mutation}
of $Q$, as follows:
(i) For any pair $a:j\to i$ and $b:i\to k$ of arrows in $Q$, we create a new arrow $[ba]:j\to k$. (ii) Reverse all arrows starting or ending at $i$.
(iii) Remove a maximal disjoint set of 2-cycles.
For example, if $Q$ is $1\xrightarrow{}2\xrightarrow{}3$, then $\mu_3(Q)$ is
$1\xrightarrow{}2\xleftarrow{}3$ and $\mu_2(Q)$ is
$\xymatrix@R=0em@C=1em{1\ar@/_4mm/[rr]&2\ar[l]&3\ar[l]}$.

\medskip
Fix a function field $F=\mathbb{Q}(x_1,\ldots,x_n)$ in $n$ variables over the field of rational
numbers $\mathbb{Q}$. Then $(\{x_1,\ldots,x_n\},Q)$ is an \emph{initial seed}, consisting of the pair
of the free generating set $\{x_1,\ldots,x_n\}$ of $F$ over $\mathbb{Q}$ and the acyclic quiver $Q$.
Then for each $i=1,\ldots,n$ we define a new \emph{seed}
\[\mu_i(\{x_1,\ldots,x_n\},Q):=(\{x_1,\ldots,x_i^*,\ldots,x_n\},\mu_i(Q)),\]
where $x_i^*=\frac{m_1+m_2}{x_i}$, with $m_1$ and $m_2$ being certain monomials in
$x_1,\ldots,x_{i-1},x_{i+1},\ldots,x_n$ determined by the quiver in the seed.
We perform this process for all $\mu_1,\ldots,\mu_n$, and then on the new seeds, etc. 
It is easy to see that $\mu_i^2={\rm id}$ holds for $i=1,\ldots,n$.
This gives rise to a graph where the vertices are the seeds,
and the edges between them are induced by the $\mu_i$.
The graph may be finite or infinite.

The $n$ element subsets of $F$ obtained in this way are called \emph{clusters}, and the elements of clusters are called \emph{cluster variables}.
The associated \emph{cluster algebra} is the subalgebra of $F$ generated by
all the cluster variables.
There is only a finite number of clusters (or equivalently cluster variables)
if and only if the quiver $Q$ is Dynkin.

One approach to the study of cluster algebras is via so-called \emph{categorification}.
This means that we want to find some nice categories, like module categories or triangulated
categories,
where we have some objects with similar properties as clusters and cluster variables.
In particular we are looking for indecomposable objects corresponding to cluster variables.

As we have seen in section 1, the tilting modules almost have the desired properties,
but the operation $\mu_i$ on tilting modules is not always defined.
This motivated the introduction of cluster categories $\CC_Q$ associated with a finite
acyclic quiver \cite{BMRRT}. These categories have more indecomposable objects than
$\mod kQ$, and also more morphisms between the old objects. They are defined as orbit
categories $\DDD^{\bo}(kQ)/\tau[-1]$, where $\DDD^{\bo}(kQ)$ is the bounded derived category of $\mod kQ$ and $\tau:\DDD^{\bo}(kQ)\to\DDD^{\bo}(kQ)$ is a derived AR-translation.
For example, the AR-quiver of the cluster category of type $A_3$ was given in section 1.
The category $\CC_Q$ is Hom-finite, and also triangulated \cite{K}, and we have a
functorial isomorphism $D\Ext^1_{\CC_Q}(X,Y)\simeq\Ext^1_{\CC_Q}(Y,X)$ \cite{BMRRT}.
In general a $k$-linear Hom-finite triangulated category $\CC$ is called \emph{2-Calabi-Yau}
if there is a functorial isomorphism $D\Ext^1_{\CC}(X,Y)\simeq\Ext^1_{\CC}(Y,X)$.
Hence $\CC_Q$ is 2-Calabi-Yau.

We can regard the tilting $kQ$-modules as special objects in $\CC_Q$,
but this class of objects is not large enough
since mutation is not always possible.
Instead, if we consider all tilting modules coming from algebras derived
equivalent to $kQ$, we actually get a class of objects 
such that mutation is always possible.
These objects can be described directly as follows.

We say that an object $T$ in a 2-Calabi-Yau triangulated category $\CC$ is \emph{rigid} if
$\Ext^1_{\CC}(T,T)=0$. It is \emph{maximal rigid} if 
it is rigid, and moreover if $T\oplus X$ is rigid for some $X$ in $\CC$,
then $X$ is in $\add T$. A rigid object $T$ in $\CC$ is \emph{cluster-tilting}
if $\Ext^1_{\CC}(T,X)=0$ for some $X$ in $\CC$ implies that $X$ is in $\add T$. 
Clearly any cluster-tilting object is maximal rigid, and the converse is true
for the cluster category \cite{BMRRT}, and also for more general 2-Calabi-Yau
triangulated categories (see Corollary \ref{application to CT}(b)).
A basic object $U$ in $\CC$ is an \emph{almost complete cluster-tilting object}
if there exists an indecomposable object $X$ in $\CC$ such that $U\oplus X$ is
a basic cluster tilting object. In this case $X$ is as before called
a \emph{complement} for $U$.

We have the following result for cluster-tilting objects in
2-Calabi-Yau triangulated categories, which improves
Theorem \ref{complements for tilting} for tilting modules.

\begin{theorem}\label{complements for cluster-tilting}\emph{\cite{BMRRT,IY}}
Let $\CC$ be a 2-Calabi-Yau triangulated category.
Any almost complete cluster-tilting object in $\CC$ has precisely two complements.
\end{theorem}

It was mentioned in section 1 that when we deal with support tilting modules over
path algebras $kQ$, then the mutation $\mu_i$ is defined for all $i=1,\ldots,n$.
To prove this, one uses the corresponding result for cluster-tilting objects
in cluster categories and interprets the condition $\Ext^1_{\CC_Q}(T,T)=0$
and $|T|=n$ for $T$ in a cluster category as a condition for $\mod kQ$.

Further basic results about cluster categories were developed.
They were of interest in themselves, and also helped to establish
a closer connection with acyclic cluster
algebras, in a series of papers including \cite{BMR,BMRT,CC,CK1,CK2}.
This led to the result that there is a map $\CC_Q\to A$ for the
acyclic cluster algebra $A$ associated with the quiver $Q$,
called \emph{cluster character}, inducing a bijection between indecomposable
rigid objects in $\CC_Q$ and the cluster variables in $A$.

We now give some classes of 2-Calabi-Yau triangulated categories with cluster-tilting objects,
in addition to the cluster categories.

(1) \cite{IR,BIRS} Let $W_Q$ be the Coxeter group associated with a finite acyclic quiver $Q$.
For any element $w\in W_Q$, there is associated a 2-Calabi-Yau triangulated
category $\CC_w$ as follows.
Let $\Pi$ be the preprojective algebra associated with $Q$
(see section 4.4).
Let $e_1,\ldots,e_n$ denote the idempotent elements associated with the vertices $1,\ldots,n$.
Then we associate with $i$ the ideal $I_i=\Pi(1-e_i)\Pi$ in $\Pi$, and the ideal
$I_w=I_{i_1}\cdots I_{i_t}$ with an element $w\in W_Q$ with a reduced expression
$w=s_{i_1}\cdots s_{i_t}$.
Let $\Pi_w=\Pi/I_w$ and $\Sub\Pi_w$ be the category of $\Pi_w$-modules which are
submodules of a finite direct sum of copies of $\Pi_w$. Then the stable category
$\CC_w:=\underline{\Sub} \Pi_w$ is a 2-Calabi-Yau triangulated category.
In particular, the stable category $\underline{\mod}\Pi$, where $\Pi$ is
the preprojective algebra of a Dynkin quiver, belongs to this class. It was studied
by Geiss, Leclerc and Schr\"oer \cite{GLS2}, who also independently studied a subclass of
$\CC_w$ \cite{GLS1}.
We will see in section 4.4 that, in the Dynkin case, the ideals $I_w$ play
a role in $\tau$-tilting theory.

(2) There is a generalization of the cluster categories $\CC_Q$ to so-called \emph{generalized
cluster categories} $\CC_\Lambda$ associated with algebras $\Lambda$ with $\gl\Lambda\le 2$
and $\CC_{Q,W}$ associated with quivers with potential $(Q,W)$ \cite{Am}.
Under a certain finiteness condition on $A$ and $(Q,W)$, the generalized cluster category
is a 2-Calabi-Yau triangulated category with cluster-tilting object.

(3) When $R$ is a 3-dimensional isolated Gorenstein singularity with $k\subset R$, then the
stable category $\underline{\CM}R$ of maximal Cohen-Macaulay $R$-modules is a 2-Calabi-Yau
triangulated category. More concrete examples are invariant rings $R=k[[X,Y,Z]]^G$,
where $G$ is a finite subgroup of ${\rm SL}_3(k)$ acting freely on $k^3\backslash\{0\}$.
For other examples, see \cite{BIKR}.
These 2-Calabi-Yau triangulated categories have quite different origin from those in (2),
but it was shown in \cite{AmIR} that they are often equivalent.

\section{4. $\tau$-tilting theory}

In this section we discuss the central results of $\tau$-tilting theory, including those
analogous to results from tilting theory. Let as before
$\Lambda$ be a finite dimensional $k$-algebra.

\subsection{4.1. Definition and basic properties}\label{section: 4.1}

\ We have dualities
$D=\Hom_k(-,k):\mod\Lambda\leftrightarrow\mod\Lambda^{\op}$ and 
$(-)^*:=\Hom_\Lambda(-,\Lambda):\proj\Lambda\leftrightarrow\proj\Lambda^{\op}$
which induce an equivalence $\nu:=D(-)^*:\proj\Lambda\to\inj\Lambda$
called \emph{Nakayama functor}.
For $X$ in $\mod\Lambda$ with a minimal projective presentation
$P_1\xrightarrow{d_1}P_0\xrightarrow{d_0}X\xrightarrow{}0$,
we define $\Tr X$ in $\mod\Lambda^{\op}$ and $\tau X$ in $\mod\Lambda$ by the exact sequences
\[P_0^*\xrightarrow{d_1^*}P_1^*\xrightarrow{}\Tr X\xrightarrow{}0\ \ \mbox{and}\ \ 
0\xrightarrow{}\tau X\xrightarrow{}\nu P_1\xrightarrow{\nu d_1}\nu P_0.\]
Then $\tau$ (respectively, $\Tr$) gives bijections between the isomorphism classes of
indecomposable non-projective $\Lambda$-modules and
the isomorphism classes of indecomposable non-injective $\Lambda$-modules
(respectively, non-projective $\Lambda^{\op}$-modules).
More strongly, $\tau$ and $\Tr$ give an equivalence and a duality respectively
between the stable categories.
We denote by $\underline{\mod}\Lambda$ the \emph{stable category} modulo projectives and by 
$\overline{\mod}\Lambda$ the \emph{costable category} modulo injectives.
Then $\Tr$ gives a duality $\Tr:\underline{\mod}\Lambda\leftrightarrow
\underline{\mod}\Lambda^{\op}$ called the \emph{transpose},
and $\tau$ gives an equivalence $\tau:\underline{\mod}\Lambda\to\overline{\mod}\Lambda$
called the \emph{AR-translation}. Moreover we have a functorial isomorphism
\[\underline{\Hom}_\Lambda(X,Y)\simeq D\Ext^1_\Lambda(Y,\tau X)\]
for any $X$ and $Y$ in $\mod\Lambda$ called \emph{AR-duality}. 

There is an important class of modules, which were introduced by Auslander-Smalo
\cite{AS} (see also \cite{ARS,ASS}) almost at the same time when tilting modules were
introduced in \cite{BB}, then later almost forgotten:
We call $T$ in $\mod\Lambda$ \emph{$\tau$-rigid} if $\Hom_{\Lambda}(T,\tau T)=0$.
In this case we have $\Ext_\Lambda^1(T,T)=0$ (i.e. $T$ is \emph{rigid}) by AR-duality.
Now let us call $T$ in $\mod\Lambda$ {\it $\tau$-tilting} if $T$ is $\tau$-rigid and $|T|=|\Lambda|$.
We call $T$ in $\mod\Lambda$ {\it support $\tau$-tilting} if there exists an idempotent $e$ of $\Lambda$ such that $T$ is a $\tau$-tilting $(\Lambda/\langle e\rangle)$-module.

Any local algebra $\Lambda$ has precisely two basic support $\tau$-tilting modules
$\Lambda$ and $0$.
It follows from AR-duality that any partial tilting module is $\tau$-rigid,
and any tilting module is $\tau$-tilting. But the converse is far from being true.
In fact, a selfinjective algebra $\Lambda$ usually has a lot of $\tau$-tilting modules
(see e.g. section 4.4) even though it has a unique basic tilting module $\Lambda$.
For another example, let $Q$ be the quiver $1\xrightarrow{a}2\xrightarrow{b}3$
and $\Lambda:=kQ/\langle ba\rangle$.
Then $T=P_1\oplus S_1\oplus P_3$ is a $\tau$-tilting $\Lambda$-module
which is not a tilting $\Lambda$-module.

The next observation gives a more precise relationship between
tilting modules and $\tau$-tilting modules, where a $\Lambda$-module
$X$ is called \emph{sincere} if any simple $\Lambda$-module appears in $X$
as a composition factor. Clearly any faithful module is sincere.

\begin{proposition}\emph{\cite{AIR}}
{\rm (a)} Tilting modules are precisely faithful (support) $\tau$-tilting modules.\\
{\rm (b)} Any $\tau$-tilting (respectively, $\tau$-rigid) $\Lambda$-module $T$ is a tilting (respectively, partial tilting) $(\Lambda/\ann T)$-module.\\
{\rm (c)} $\tau$-tilting modules are precisely sincere support $\tau$-tilting modules.
\end{proposition}

When $T$ is a $(\Lambda/\langle e\rangle)$-module for some idempotent $e$,
it is often useful to compare properties of $T$ as $(\Lambda/\langle e\rangle)$-module
with properties of $T$ as a $\Lambda$-module.

\begin{proposition}\emph{\cite{AIR}}\label{idempotent quotient}
Let $e$ be an idempotent of $\Lambda$ and $T\in\mod(\Lambda/\langle e\rangle)$.
Then $T$ is a $\tau$-rigid $\Lambda$-module if and only if
$T$ is a $\tau$-rigid $(\Lambda/\langle e\rangle)$-module.
\end{proposition}

There is a nice relationship between torsion classes and support $\tau$-tilting modules,
which improves Theorem \ref{tilting and torsion2}.
We denote by $\sttilt\Lambda$ the set of isomorphism classes of basic support
$\tau$-tilting $\Lambda$-modules, and by $\ftors\Lambda$ the set of functorially
finite torsion classes in $\mod\Lambda$.

\begin{theorem}\emph{\cite{AIR}}\label{basic bijection}
There is a bijection $\sttilt\Lambda\rightarrow\ftors\Lambda$ given by $T\mapsto\Fac T$.
\end{theorem}

This gives a natural way of extending a $\tau$-rigid module to a $\tau$-tilting
module, similar to Bongartz's result (Theorem \ref{Bongartz}):
For any $\tau$-rigid module $U$,
we have an associated functorially finite torsion class ${}^{\perp}(\tau U)$.
By Theorem \ref{basic bijection}, there exists a support $\tau$-tilting module
$T$ such that $\Fac T={}^{\perp}(\tau U)$. It is easy to show that $T$ is a
$\tau$-tilting module satisfying $U\in\add T$. Consequently we have the following result.

\begin{theorem}\emph{\cite{AIR}}\label{Bongartz for tau-rigid}
Any $\tau$-rigid module is a direct summand of a $\tau$-tilting module.
\end{theorem}

The next result gives several characterizations of $\tau$-tilting modules. This
is important because it has a result from \cite{ZZ} listed in
Corollary \ref{application to CT} as a special case.

\begin{theorem}\label{3 conditions}\emph{\cite{AIR}}
The following conditions are equivalent for a $\tau$-rigid $\Lambda$-module $T$.\\
{\rm (a)} $T$ is $\tau$-tilting.\\
{\rm (b)} $T$ is maximal $\tau$-rigid, i.e. if $T\oplus X$ is $\tau$-rigid for some $\Lambda$-module $X$, then $X\in\add T$.\\
{\rm (c)} $^{\perp}(\tau T)=\Fac T$.\\
{\rm (d)} If $\Hom_\Lambda(T,\tau X)=0$ and $\Hom_\Lambda(X,\tau T)=0$, then $X\in\add T$.
\end{theorem}

In addition to $\tau$-rigid and support $\tau$-tilting modules, it is also convenient
to consider associated pairs of modules.

Let $(T,P)$ be a pair with $T\in\mod\Lambda$ and $P\in\proj\Lambda$.
We call $(T,P)$ a {\it $\tau$-rigid pair} if $T$ is $\tau$-rigid and $\Hom_\Lambda(P,T)=0$.
We call $(T,P)$ a {\it support $\tau$-tilting pair} if
$(T,P)$ is $\tau$-rigid and $|T|+|P|=|\Lambda|$.
The notion of support $\tau$-tilting module and that of support $\tau$-tilting pair
are essentially the same since, for an idempotent $e$ of $\Lambda$, the pair
$(T,\Lambda e)$ is support $\tau$-tilting if and only if $T$ is a $\tau$-tilting
$(\Lambda/\langle e\rangle)$-module (cf. Proposition \ref{idempotent quotient}).
Moreover $e$ is uniquely determined by $T$ in this case.

Next we discuss the left-right symmetry of support $\tau$-tilting modules.
This is somewhat surprising since it does not have an analog for tilting modules.
We decompose $T$ in $\mod\Lambda$ as $T=T_{\rm pr}\oplus T_{\rm np}$ where $T_{\rm pr}$ is a maximal projective direct summand of $T$.
For a $\tau$-rigid pair $(T,P)$ for $\Lambda$, let
\[(T,P)^\dagger:=((\Tr T_{\rm np})\oplus P^*,T_{\rm pr}^*).\]
We denote by $\trigid\Lambda$ the set of isomorphism classes of basic $\tau$-rigid pairs of 
$\Lambda$.

\begin{theorem}\label{left-right symmetry}\emph{\cite{AIR}}
$(-)^\dagger$ gives bijections
$\trigid\Lambda\leftrightarrow\trigid\Lambda^{\op}$
and $\sttilt\Lambda\leftrightarrow\sttilt\Lambda^{\op}$
such that $(-)^{\dagger\dagger}={\rm id}$.
\end{theorem}

\subsection{4.2. Mutation and Partial order}

\ Now we give a basic theorem in the mutation theory for support $\tau$-tilting modules,
which improves Theorem \ref{complements for tilting} for tilting modules.
We say that a $\tau$-rigid pair $(U,P)$ is \emph{almost complete support $\tau$-tilting}
if $|U|+|P|=|\Lambda|-1$.
In this case, for an indecomposable module $X$, we say that $(X,0)$ (respectively, $(0,X)$)
is a \emph{complement} of $(U,P)$ if $(U\oplus X,P)$ (respectively, $(U,P\oplus X)$)
is a support $\tau$-tilting pair, which we call a \emph{completion} of $(U,P)$.

\begin{theorem}\label{2 complements}\emph{\cite{AIR}}
Let $\Lambda$ be a finite dimensional $k$-algebra.
Then any almost complete support $\tau$-tilting pair has precisely two complements.
\end{theorem}

Two completions $(T,Q)$ and $(T',Q')$ of an almost complete support $\tau$-tilting pair
$(U,P)$ are called \emph{mutations} of each other.
We write $(T',Q')=\mu_{(X,0)}(T,Q)$ (respectively, $(T',Q')=\mu_{(0,X)}(T,Q)$) if
$(X,0)$ (respectively, $(0,X)$) is a complement of $(U,P)$ giving rise to $(T,Q)$.

For example, let $Q$ be the quiver $1\xrightarrow{a}2\xrightarrow{b}3$
and $\Lambda:=kQ/\langle ba\rangle$.
Then $\mu_{(P_1,0)}(\Lambda,0)=(P_2\oplus P_3,P_1)$,
$\mu_{(P_2,0)}(\Lambda,0)=(P_1\oplus S_1\oplus P_3,0)$, 
and $\mu_{(P_3,0)}(\Lambda,0)=(P_1\oplus P_2\oplus S_2,0)$.

Note that it is easy to understand when we have a complement of $(U,P)$ of the form $(0,X)$
and how to find it.
This occurs precisely when $U$ is not a sincere $(\Lambda/\langle e\rangle)$-module
for an idempotent $e$ of $\Lambda$ satisfying $\add P=\add(\Lambda e)$.
In this case there exists a primitive idempotent $e'$ of $\Lambda$ such that
$U$ is a sincere $(\Lambda/\langle e+e'\rangle)$-module, and $(0,\Lambda e')$
is a complement of $(U,P)$.

The next result gives a relationship between complements of the form $(X,0)$,
which is an analog of Theorem \ref{exchange sequence for tilting} for
almost complete tilting modules.

\begin{theorem}\label{exchange sequence for tau-tilting}\emph{\cite{AIR}}
Let $(X,0)$ and $(Y,0)$ be two complements of an almost complete support $\tau$-tilting pair
$(U,P)$.
After we interchange $X$ and $Y$ if necessary, there is an exact sequence
$X\xrightarrow{f}U'\xrightarrow{g}Y^r\to0$ for some positive integer $r$, where
$f$ is a minimal left $(\add U)$-approximation and $g$ is a minimal right
$(\add U)$-approximation.
\end{theorem}

A main difference from Theorem 3 is that $f$ above is not necessarily injective.
We conjecture that the integer $r$ appearing in the exact sequence above is always one.

Now we introduce an important partial order on the set $\sttilt\Lambda$.
We have a bijection $\sttilt\Lambda\to\ftors\Lambda$ given by
$T\mapsto\Fac T$ in Theorem \ref{basic bijection}. Then inclusion in $\ftors\Lambda$
gives rise to a partial order on $\sttilt\Lambda$, i.e. we define
\[T\ge U\]
if and only if $\Fac T\supset\Fac U$.
Clearly $(\Lambda,0)$ is the maximum element and $(0,\Lambda)$ is the minimum element 
of $\sttilt\Lambda$. We give basic properties.

\begin{proposition}\label{basic of partial order}\emph{\cite{AIR}}\\
{\rm (a)} The bijection $(-)^\dagger:\sttilt\Lambda\to\sttilt\Lambda^{\rm op}$ reverses
the partial order.\\
{\rm (b)} Let $(T,Q)$ and $(T',Q')$ be support $\tau$-tilting pairs for 
$\Lambda$ which are mutations of each other.
Then we have either $(T,Q)>(T',Q')$ or $(T,Q)<(T',Q')$.
\end{proposition}

Now we have a recipe to calculate mutation of support $\tau$-tilting pairs,
which is slightly more complicated than that for tilting modules because the map $f$
in Theorem \ref{exchange sequence for tau-tilting} is not necessarily injective.
Let $(U,P)$ be an almost complete support $\tau$-tilting pair with a completion $(T,Q)$.
We give a method to calculate another completion $(T',Q')$ of $(U,P)$ from $(T,Q)$,
which is divided into two cases (A) and (B) below, depending on which of $(T,Q)$ and
$(T',Q')$ is bigger. This can be checked without knowing $(T',Q')$ since $(T,Q)>(T',Q')$ holds
if and only if $T=U\oplus X$ for an indecomposable module $X$ satisfying $X\notin\Fac U$.\\
{\rm (A)} Assume $(T,Q)>(T',Q')$ holds. Then we take a minimal left $(\add U)$-approximation
$X\xrightarrow{f}U'$ of $X$. If $f$ is not surjective, then the cokernel of $f$ is a
direct sum of copies of an indecomposable module $Y$
(Theorem \ref{exchange sequence for tau-tilting}), and we have $(T',Q')=(U\oplus Y,P)$.
If $f$ is surjective, then there exists an indecomposable projective module $Y$
such that $(T',Q')=(U,P\oplus Y)$.\\
{\rm (B)} Assume $(T,Q)<(T',Q')$ holds. Since $(T,Q)^\dagger>(T',Q')^\dagger$
holds by Proposition \ref{basic of partial order}(a), we can calculate
$(T',Q')^\dagger$ from $(U,P)^\dagger$ by using the left approximation, as we
observed in (A). Applying $(-)^\dagger$ to $(T',Q')^\dagger$, we obtain $(T',Q')$.

For example, let $Q$ be the quiver $1\xrightarrow{a}2\xrightarrow{b}3$
and $\Lambda:=kQ/\langle ba\rangle$.
First we caclulate $\mu_{(P_2,0)}(\Lambda,0)$. 
Since $P_2\notin\Fac(\Lambda/P_2)$, we use (A) above.
The minimal left $(\add\Lambda/P_2)$-approximation
of $P_2$ is $P_2\to P_1$, which has cokernel $S_1$.
Thus $\mu_{(P_2,0)}(\Lambda,0)=(P_1\oplus S_1\oplus P_3,0)$.

Next we calculate $\mu_{(0,P_3)}(P_1\oplus S_1,P_3)$.
We must use (B) in this case.
We have $(P_1\oplus S_1,P_3)^\dagger=((\Tr S_1)\oplus P_3^*,P_1^*)=
(S'_2\oplus P_3^*,P_1^*)$, where $S'_i$ is the simple $\Lambda^{\op}$-module
associated with the vertex $i$.
Now we have $\mu_{(P_3^*,0)}(S'_2\oplus P_3^*,P_1^*)=
(S'_2,P_1^*\oplus P_3^*)$.
Thus $\mu_{(0,P_3)}(P_1\oplus S_1,P_3)=
(S'_2,P_1^*\oplus P_3^*)^\dagger=(P_1\oplus S_1\oplus P_3,0)$.

The following result shows a strong connection between mutation and partial order.

\begin{theorem}\label{Hasse is mutation}\emph{\cite{AIR}}
For $T,U\in\sttilt\Lambda$, the following conditions are equivalent.\\
{\rm (a)} $T$ and $U$ are mutations of each other, and $T>U$.\\
{\rm (b)} $T>U$ and there is no $V\in\sttilt\Lambda$ such that $T>V>U$.
\end{theorem}

We define the \emph{exchange quiver} of $\sttilt\Lambda$ as the quiver whose
set of vertices is $\sttilt\Lambda$ and we draw an arrow from $T$ to $U$ if
$U$ is a mutation of $T$ such that $T>U$.
The following analog of Theorem \ref{Hasse and exchange}
is an easy consequence of Theorem \ref{Hasse is mutation}.

\begin{corollary}\label{Hasse is mutation 2}\emph{\cite{AIR}}
The exchange quiver of $\sttilt\Lambda$ coincides with
the Hasse quiver of the partially ordered set $\sttilt\Lambda$.
\end{corollary}

Another application of Theorem \ref{Hasse is mutation} is the following.

\begin{corollary}\label{connected}\emph{\cite{AIR}}
If the support $\tau$-tilting quiver has a finite
connected component, then it is connected.
\end{corollary}


We end this subsection with the following generalization of Theorem
\ref{complements for cluster-tilting}, where $\Gamma$ is given
naturally by using Bongartz completion in Theorem \ref{Bongartz for tau-rigid}.

\begin{theorem}\emph{\cite{J}}\label{jasso}
Let $(U,P)$ be a $\tau$-rigid pair for $\Lambda$, and
$\sttilt_{(U,P)}\Lambda$ the subset of $\sttilt\Lambda$
consisting of support $\tau$-tilting pairs which have $(U,P)$ as a direct summand.
Then there is a bijection $\sttilt\Gamma\to\sttilt_{(U,P)}\Lambda$
for some finite dimensional $k$-algebra $\Gamma$ with
$|\Gamma|=|\Lambda|-|U|-|P|$.
\end{theorem}

When $(U,P)$ is an almost complete support $\tau$-tilting pair, the algebra
$\Gamma$ is local and hence $\sttilt\Gamma=\{\Gamma,0\}$ holds. Thus Theorem
\ref{complements for cluster-tilting} is a special case of Theorem \ref{jasso}.

\subsection{4.3. Connection with silting complexes}\label{section: silting}

\ In this subsection, we observe that there are bijections between
support $\tau$-tilting modules and a certain class of silting complexes
defined in section 2.2. Let $\Lambda$ be a finite dimensional $k$-algebra.
We call a complex $P=(P^{i},d^i)$ in $\KKb(\proj \Lambda)$
\emph{two-term} if $P^{i}=0$ for all $i\neq 0,-1$.
We denote by $\twosilt\Lambda$ the set of isomorphism classes of basic two-term
silting complexes in $\KKb(\proj \Lambda)$.
Two-term tilting complexes were studied by Hoshino, Kato and Miyachi \cite{HKM}.
We have the following connection between support $\tau$-tilting modules and two-term
silting complexes.

\begin{theorem}\label{sptsil4}\emph{\cite{AIR,DF}}
Let $\Lambda$ be a finite dimensional $k$-algebra.
Then there exists a bijection $\twosilt\Lambda\rightarrow\sttilt\Lambda$ given by
$P\mapsto H^0(P)$.
\end{theorem}

The natural bijection $\twosilt\Lambda\rightarrow\twosilt\Lambda^{\op}$ given by
$P\mapsto\Hom_\Lambda(P,\Lambda)$ corresponds to the bijection in
Theorem \ref{left-right symmetry}.

Immediately we have the following consequence.

\begin{corollary}\emph{\cite{AIR,DF}}\label{complement for two-term silting}
Let $\Lambda$ be a finite dimensional $k$-algebra.
Then any two-term almost complete silting complex has precisely
two complements which are two-term.
\end{corollary}

The following property of $\tau$-rigid pairs is an analog of a property of silting complexes \cite[2.25]{AI}.

\begin{proposition}\emph{\cite{AIR}}
Let $(T,Q)$ be a $\tau$-rigid pair for $\Lambda$ and
$Q_1\to Q_0\to T\to0$ a minimal projective presentation.
Then $Q_0$ and $Q_1\oplus Q$ have no non-zero direct summands in common.
\end{proposition}

Let $K_0(\proj\Lambda)$ be the Grothendieck group of the additive category
$\proj\Lambda$ of finitely generated projective $\Lambda$-modules.
It is a free abelian group with a basis consisting of the isomorphism classes
$P_1,\ldots,P_n$ of indecomposable projective $\Lambda$-modules.
Let $Q_1\to Q_0\to T\to0$ be a minimal projective presentation of $T$ in $\mod\Lambda$.
The element
\[g^T:=Q_0-Q_1\]
in $K_0(\proj\Lambda)$ is called the \emph{index} 
or the \emph{$g$-vector} of $T$.
The following result is an immediate consequence of Theorem \ref{sptsil4} and
a general result for silting complexes \cite[2.27]{AI}.

\begin{theorem}\label{g-vectors}\emph{\cite{AIR}}
Let $(T_1\oplus\cdots\oplus T_{\ell},Q_{\ell+1}\oplus\cdots\oplus Q_n)$
be a support $\tau$-tilting pair for $\Lambda$ with $T_i$ and $Q_i$ indecomposable.
Then $g^{T_1},\cdots,g^{T_\ell},g^{Q_{\ell+1}},\cdots,g^{Q_n}$ form a basis of
$K_0(\proj\Lambda)$.
\end{theorem}

The following result shows that $g$-vectors determine $\tau$-rigid pairs,
which is an analog of \cite[2.3]{DK}.

\begin{theorem}\label{g-vectors are complete invariant}\emph{\cite{AIR}}
We have an injective map from the set of isomorphism classes of $\tau$-rigid
pairs for $\Lambda$ to $K_0(\proj\Lambda)$ given by $(T,P)\mapsto g^T-g^P$.
\end{theorem}

\subsection{4.4. Examples}

\ In this section, we use a convention to describe modules via their composition
series. For example, $\begin{smallmatrix}1\\ 2\end{smallmatrix}$ is an indecomposable
$\Lambda$-module $X$ with a simple submodule $S_2$ such that $X/S_2=S_1$.
For example, each simple module $S_i$ is written as $\begin{smallmatrix}i\end{smallmatrix}$.

Let $Q$ be the quiver $1\xrightarrow{a}2\xrightarrow{b}3$ and $\Lambda=kQ/\langle ba\rangle$.
Then $P_1=\begin{smallmatrix}1\\ 2\end{smallmatrix}$,
$P_2=\begin{smallmatrix}2\\ 3\end{smallmatrix}$ and
$P_3=\begin{smallmatrix}3\end{smallmatrix}$ in this case.
The Hasse quiver of $\sttilt\Lambda$ is the following.
{\small 
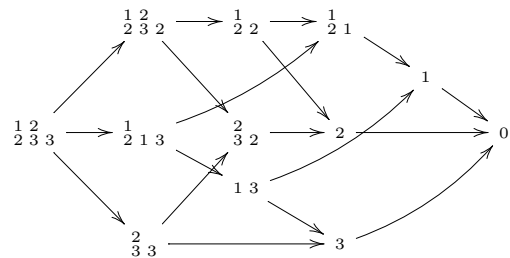
\[\xymatrix@R1em@C2em{
&{\begin{smallmatrix}1&2&\\ 2&3&2\end{smallmatrix}}\ar[r]\ar[rdd]&
{\begin{smallmatrix}1&\\ 2&2\end{smallmatrix}}\ar[r]\ar[rdd]&
{\begin{smallmatrix}1&\\ 2&1\end{smallmatrix}}\ar[rd]\\
&&&&{\begin{smallmatrix}1\end{smallmatrix}}\ar[rd]\\
{\begin{smallmatrix}1&2&\\ 2&3&3\end{smallmatrix}}\ar[r]\ar[ruu]\ar[rdd]&
{\begin{smallmatrix}1&&\\ 2&1&3\end{smallmatrix}}\ar@/_2mm/[rruu]\ar[rd]&
{\begin{smallmatrix}2&\\ 3&2\end{smallmatrix}}\ar[r]&
{\begin{smallmatrix}2\end{smallmatrix}}\ar[rr]&&
{\begin{smallmatrix}0\end{smallmatrix}}\\
&&{\begin{smallmatrix}1&3\end{smallmatrix}}\ar@/_2mm/[rruu]\ar[rd]\\
&{\begin{smallmatrix}2&\\ 3&3\end{smallmatrix}}\ar[ruu]\ar[rr]&&
{\begin{smallmatrix}3\end{smallmatrix}}\ar@/_2mm/[rruu]
}\]}
where the symbol $\oplus$ is omitted.
For example, $\begin{smallmatrix}1&2&\\ 2&3&3\end{smallmatrix}$ is
$\begin{smallmatrix}1\\ 2\end{smallmatrix}\oplus\begin{smallmatrix}2\\ 3\end{smallmatrix}
\oplus\begin{smallmatrix}3\end{smallmatrix}=P_1\oplus P_2\oplus P_3$.

Let $Q$ be the quiver $1\xto{a}2\xto{a}3\xto{a}1$ and
$\Lambda=kQ/\langle a^2=0\rangle$.
Then the Hasse quiver of $\sttilt\Lambda$ is the following.
{\small
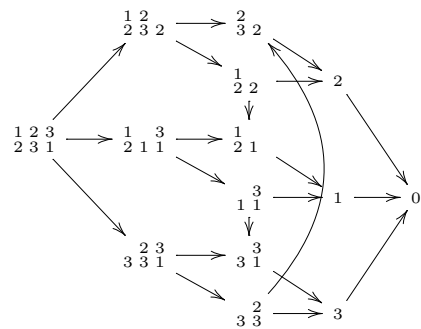
\[\xymatrix@R.8em@C2em{
&{\begin{smallmatrix}1&2&\\ 2&3&2\end{smallmatrix}}\ar[r]\ar[rd]&
{\begin{smallmatrix}2&\\ 3&2\end{smallmatrix}}\ar[rd]\\
&&{\begin{smallmatrix}1&&\\ 2&2\end{smallmatrix}}\ar[r]\ar[d]&
{\begin{smallmatrix}2\end{smallmatrix}}\ar[rdd]\\
{\begin{smallmatrix}1&2&3\\ 2&3&1\end{smallmatrix}}\ar[ruu]\ar[r]\ar[rdd]&
{\begin{smallmatrix}1&&3\\ 2&1&1\end{smallmatrix}}\ar[r]\ar[rd]&
{\begin{smallmatrix}1&&\\ 2&1\end{smallmatrix}}\ar[rd]\\
&&{\begin{smallmatrix}&3\\ 1&1\end{smallmatrix}}\ar[r]\ar[d]&
{\begin{smallmatrix}1\end{smallmatrix}}\ar[r]&
{\begin{smallmatrix}0\end{smallmatrix}}\\
&{\begin{smallmatrix}&2&3\\ 3&3&1\end{smallmatrix}}\ar[r]\ar[rd]&
{\begin{smallmatrix}&3\\ 3&1\end{smallmatrix}}\ar[rd]\\
&&{\begin{smallmatrix}&2\\ 3&3\end{smallmatrix}}\ar[r]\ar@/_10mm/[uuuuu]&
{\begin{smallmatrix}3\end{smallmatrix}}\ar[ruu]
}\]}
This coincides with the exchange graph of cluster-tilting objects
in the cluster category of type $A_3$ given in section 1,
which is a consequence of Theorem \ref{bijection between CT and PT} below.

Now we deal with an important class of algebras.
Let $Q$ be a Dynkin quiver.
Define a new quiver $\overline{Q}$ by adding a new arrow $a^*:j\to i$ to $Q$
for each arrow $a:i\to j$ in $Q$.
We call
\[\Pi:=K\overline{Q}/\langle\sum_{a:\mbox{{\scriptsize arrow in }}Q}(aa^*-a^*a)\rangle\]
the \emph{preprojective algebra} of $Q$.

Any preprojective algebra $\Pi$ of Dynkin type is a finite dimensional
algebra which is selfinjective. In particular $\Pi$ is a unique
tilting $\Pi$-module. On the other hand we will see that
there is a large family of support $\tau$-tilting $\Pi$-modules.

For a Dynkin quiver $Q$, we denote by $W_Q$ the associated \emph{Weyl group}, i.e.
$W_Q$ is presented by generators $s_1,\ldots,s_n$ with the following relations:\\
$\bullet$ $s_i^2=1$,\\
$\bullet$ $s_is_j=s_js_i$ if there is no arrow between $i$ and $j$ in $Q$,\\
$\bullet$ $s_is_js_i=s_js_is_j$ if there is precisely one arrow between $i$ and $j$ in $Q$.

We say that an expression $w=s_{i_1}\cdots s_{i_\ell}$ of $w\in W_Q$
with $i_1,\ldots,i_\ell\in\{1,\ldots,n\}$ is \emph{reduced} if $\ell$
is the smallest possible number. Such $\ell$ is called the \emph{length}
of $w$ and denoted by $\ell(w)$.
We have a partial order on $W_Q$, called the \emph{right order}:
We define $w\ge w'$ if and only if $\ell(w)=\ell(w')+\ell(w'^{-1}w)$.

When $Q$ is a Dynkin quiver, there is the following nice connection between elements
of $W_Q$ and support $\tau$-tilting $\Pi$-modules.
A similar result holds for tilting modules over preprojective algebras
of non-Dynkin type \cite{IR,BIRS}.

\begin{theorem}\emph{\cite{Miz}}
Let $\Pi$ be a preprojective algebra of Dynkin type
and $I_i$ is an ideal $\Pi(1-e_i)\Pi$ in $\Pi$.\\
{\rm (a)} There exists a bijection $W_Q\to\sttilt\Pi$
given by $w=s_{i_1}\cdots s_{i_\ell}\mapsto I_w:=I_{i_1}\cdots I_{i_\ell}$,
where $w=s_{i_1}\cdots s_{i_\ell}$ is any reduced expression of $w$.\\
{\rm (b)} The right order on $W_Q$ coincides with the partial order on $\sttilt\Pi$.
\end{theorem}

The Hasse quiver of $\sttilt\Pi$ for type $A_2$ is the following.
\[\xymatrix@C2em@R0em{
&{\begin{smallmatrix}1\\ 2&1\end{smallmatrix}}\ar[r]&
{\begin{smallmatrix}1\end{smallmatrix}}\ar[rd]&\\
{\begin{smallmatrix}1&2\\ 2&1\end{smallmatrix}}\ar[ur]\ar[dr]&&&
{\begin{smallmatrix}0\end{smallmatrix}}\\
&{\begin{smallmatrix}&2\\ 2&1\\ \end{smallmatrix}}\ar[r]&
{\begin{smallmatrix}2\end{smallmatrix}}\ar[ru]&}
\]
The Hasse quiver of $\sttilt\Pi$ for type $A_2$ is the following.
\[\xymatrix@C.5em@R0em{
& & & & & {\begin{smallmatrix}& \\1 &2 \\2 &1 \end{smallmatrix}}
\ar[rrd]\ar[rrddd] & & & & & \\
& & & {\begin{smallmatrix}1 && \\2 &1 &2 \\3&2 &1
\end{smallmatrix}}\ar[rru]\ar[rrd] & & & & {\begin{smallmatrix} 1
& \\2 &1 \end{smallmatrix}}\ar[rrdd] & & & \\
& & & & & {\begin{smallmatrix} 1&& \\2 &1 & \\3&2 & 1 \end{smallmatrix}}
\ar[rru] & & & & & \\
& {\begin{smallmatrix}1 &2& \\ 2&31 &2 \\3&2 &1 \end{smallmatrix}}
\ar[rruu]\ar[rrdd] & & {\begin{smallmatrix} 1&& \\2 &31 & \\3&2 &1
\end{smallmatrix}}\ar[rru]\ar[rrddd] & & & & {\begin{smallmatrix} &2
\\ 2 &1 \end{smallmatrix}}\ar[rrdd] & & {\begin{smallmatrix} 1
\end{smallmatrix}}\ar[rdd] & \\
& & & & & {\begin{smallmatrix} && \\2 & &2 \\3& 2&1
\end{smallmatrix}}\ar[rru]\ar[rrddd] & & & & & \\
{\begin{smallmatrix} 1&2&3 \\ 2&31 &2 \\3&2 &1
\end{smallmatrix}}\ar[uur]\ar[rdd]\ar[r] & {\begin{smallmatrix} 1&&3 \\
2&31 & 2\\3&2 &1 \end{smallmatrix}}\ar[rruu]\ar[rrdd] &
&{\begin{smallmatrix} &2& \\ 2&13 &2 \\3&2 &1 \end{smallmatrix}}\ar[rru]
& & & & {\begin{smallmatrix} 3&1
\end{smallmatrix}}\ar[rruu]\ar[rrdd] & & {\begin{smallmatrix} 2
\end{smallmatrix}}\ar[r] & {\begin{smallmatrix} 0
\end{smallmatrix}} \\
& & & & & {\begin{smallmatrix} && \\ & 31& \\3&2 &1
\end{smallmatrix}}\ar[rru] & & & & & \\
&{\begin{smallmatrix} &2& 3\\ 2&13 &2 \\3&2 &1
\end{smallmatrix}}\ar[rruu]\ar[rrdd] & & {\begin{smallmatrix} && 3\\ &31
&2 \\3&2 &1 \end{smallmatrix}}\ar[rru]\ar[rrddd] & & & &
{\begin{smallmatrix} 2& \\ 3&2 \end{smallmatrix}}\ar[rruu] & &
{\begin{smallmatrix} 3 \end{smallmatrix}}\ar[ruu] & \\
& & & & &{\begin{smallmatrix} 2&3 \\3&2 \end{smallmatrix}}
\ar[rru]\ar[rrd] & & & & & \\
& & & {\begin{smallmatrix} &&3 \\ 2& 3&2 \\3&2 &1 \end{smallmatrix}}
\ar[rru]\ar[rrd] & & & &{\begin{smallmatrix} & 3 \\ 3&2
\end{smallmatrix}}\ar[rruu] & & & \\
& & & & & {\begin{smallmatrix} &&3 \\ &3 & 2\\3&2 &1
\end{smallmatrix}}\ar[rru] & & & & &
}\]
These quivers coincide with the Hasse quivers of the symmetric groups
${\mathfrak S}_3$ and ${\mathfrak S}_4$ respectively.

Finally let $Q$ be the cyclic quiver $1\xto{a}2\xto{a}\cdots\xto{a}n\xto{a}1$
and $\Lambda=kQ/\langle a^m\rangle$, where $n$ and $m$ are arbtrary
positive integer satisfying $m\ge n$.
Then we have the following result, where $\ttilt\Lambda$ is the set
of isomorphism classes of basic $\tau$-tilting $\Lambda$-modules.

\begin{theorem}\emph{\cite{Ad}}
There are bijections between the following sets.\\
{\rm (a)} $\ttilt\Lambda$.\\
{\rm (b)} $\sttilt\Lambda-\ttilt\Lambda$.\\
{\rm (c)} The set of triangulations of an $n$-gon with a puncture.\\
{\rm (d)} The set of sequences $(a_1,\ldots,a_n)$ of non-negative integers
satisfying $\sum_{i=1}^na_i=n$.\\
\end{theorem}

The Hasse quiver of $\sttilt\Lambda$ for $m=n=3$ is the following.
{\small
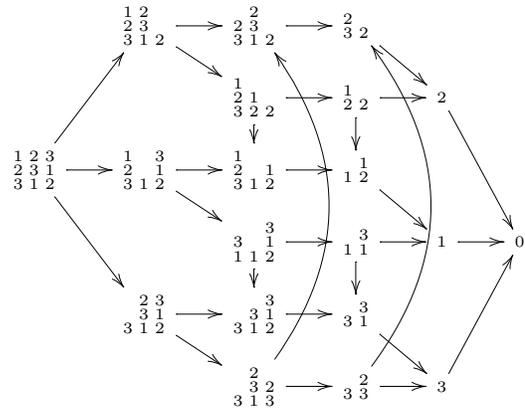
\[\xymatrix@R.8em@C2em{
&{\begin{smallmatrix}1&2&\\ 2&3&\\ 3&1&2\end{smallmatrix}}\ar[r]\ar[rd]&
{\begin{smallmatrix}&2&\\ 2&3&\\ 3&1&2\end{smallmatrix}}\ar[r]&
{\begin{smallmatrix}2&\\ 3&2\end{smallmatrix}}\ar[rd]\\
&&{\begin{smallmatrix}1&&\\ 2&1&\\ 3&2&2\end{smallmatrix}}\ar[r]\ar[d]&
{\begin{smallmatrix}1&\\ 2&2\end{smallmatrix}}\ar[r]\ar[d]&
{\begin{smallmatrix}2\end{smallmatrix}}\ar[rdd]\\
{\begin{smallmatrix}1&2&3\\ 2&3&1\\ 3&1&2\end{smallmatrix}}\ar[ruu]\ar[r]\ar[rdd]&
{\begin{smallmatrix}1&&3\\ 2&&1\\ 3&1&2\end{smallmatrix}}\ar[r]\ar[rd]&
{\begin{smallmatrix}1&&\\ 2&&1\\ 3&1&2\end{smallmatrix}}\ar[r]&
{\begin{smallmatrix}&1\\ 1&2\end{smallmatrix}}\ar[rd]\\
&&{\begin{smallmatrix}&&3\\ 3&&1\\ 1&1&2\end{smallmatrix}}\ar[r]\ar[d]&
{\begin{smallmatrix}&3\\ 1&1\end{smallmatrix}}\ar[r]\ar[d]&
{\begin{smallmatrix}1\end{smallmatrix}}\ar[r]&
{\begin{smallmatrix}0\end{smallmatrix}}\\
&{\begin{smallmatrix}&2&3\\ &3&1\\ 3&1&2\end{smallmatrix}}\ar[r]\ar[rd]&
{\begin{smallmatrix}&&3\\ &3&1\\ 3&1&2\end{smallmatrix}}\ar[r]&
{\begin{smallmatrix}&3\\ 3&1\end{smallmatrix}}\ar[rd]\\
&&{\begin{smallmatrix}&2&\\ &3&2\\ 3&1&3\end{smallmatrix}}\ar[r]\ar@/_10mm/[uuuuu]&
{\begin{smallmatrix}&2\\ 3&3\end{smallmatrix}}\ar[r]\ar@/_10mm/[uuuuu]&
{\begin{smallmatrix}3\end{smallmatrix}}\ar[ruu]
}\]}

\section{5. Connection with cluster-tilting theory}

In section 2 we have seen how the result that almost complete
cluster-tilting objects in cluster categories have exactly two complements
implied a similar result for support tilting modules
over path algebras. However we have seen that this does not work for any finite dimensional algebra.
The motivation for considering support $\tau$-tilting modules comes from the investigation of another
class of algebras associated with the cluster category $\CC_Q$, or more generally,
with any 2-Calabi-Yau triangulated category $\CC$ with cluster-tilting object.

Let $\CC$ be a 2-Calabi-Yau triangulated category and
$T$ a cluster-tilting object in $\CC$.
Then $\Lambda=\End_{\CC}(T)^{\op}$ is by definition
a \emph{2-Calabi-Yau tilted algebra} (a \emph{cluster-tilted algebra}
if $\CC$ is a cluster category $\CC_Q$). We then have the following.

\begin{theorem}\emph{\cite{BMR,KR}}
The functor $\overline{(-)}:=\Hom_{\CC}(T,-)$ induces
an equivalence of categories $\CC/[T[1]]\to \mod\Lambda$,
where $\CC/[T[1]]$ is the factor category of $\CC$ by the ideal $[T[1]]$
consisting of morphisms factoring through $\add T[1]$.
\end{theorem}

In this setting we can express the condition 
for an object $X$ in $\CC$ to be $\tau$-rigid in terms of the $\Lambda$-module
$\Hom_{\CC}(T,X)$. This leads us to the next result.
For an object $X$ in $\CC$,
we decompose $X=X'\oplus X''$ where $X''$ is a maximal direct summand of $X$
which belongs to $\add T[1]$.
We denote by $\ctilt\CC$ (respectively, $\rigid\CC$, $\mrigid\CC$) the set of
isomorphism classes of basic cluster-tilting (respectively, rigid, maximal rigid)
objects in $\CC$.

\begin{theorem}\label{bijection between CT and PT}\emph{\cite{AIR}}
The correspondence $X\mapsto(\overline{X'},\overline{X''[-1]})$
gives bijections $\rigid\CC\rightarrow\trigid\Lambda$ and
$\ctilt\CC\rightarrow\sttilt\Lambda$.
Moreover we have $\ctilt\CC=\mrigid\CC=\{U\in\rigid\CC\ |\ |U|=|T|\}$.
\end{theorem}

As a consequence we get as a special case the following important results.

\begin{corollary}\label{application to CT}
Let $\CC$ be a 2-Calabi-Yau triangulated category with a cluster-tilting object $T$.\\
{\rm (a)} \emph{\cite{IY}}
Any basic almost complete cluster-tilting object has precisely two complements.\\
{\rm (b)} \emph{\cite{ZZ}}
An object $U$ in $\CC$ is cluster-tilting if and only if it is maximal rigid if and only if it is rigid and $|U|=|T|$.
\end{corollary}

We also include an application to a result from \cite{DK}.
For a basic cluster-tilting object $T=T_1\oplus\cdots\oplus T_n$ in $\CC$,
we denote by $K_0(\add T)$ the Grothendieck group of the additive category $\add T$,
which is a free abelian group with a basis $T_1,\ldots,T_n$.
For an object $X$ in $\CC$, there exists a triangle $T''\to T'\to X\to T''[1]$
in $\CC$ with $T',T''$ in $\add T$. In this case we define the \emph{index}
(or \emph{$g$-vector}) of $X$ as $g_T^X:=T'-T''$ in $K_0(\add T)$.
Specializing Theorem \ref{g-vectors} to $\Lambda$ being 2-Calabi-Yau tilted,
we can translate to
considering cluster-tilting objects in 2-Calabi-Yau triangulated categories
to get the following.

\begin{corollary}\emph{\cite{DK}}
Let $\CC$ be a 2-Calabi-Yau triangulated category, and $T$ and 
$U=U_1\oplus\cdots\oplus U_n$ be
basic cluster-tilting objects with $U_i$ indecomposable.
Then $g_T^{U_1},\ldots,g_T^{U_n}$ form a basis of $K_0(\add T)$.
\end{corollary}

\end{article}








\end{document}